\documentclass[a4paper,11pt]{article}
\usepackage{titling}
\usepackage{hyperref}
\hypersetup{colorlinks=true,linkcolor=blue,urlcolor=blue,citecolor=blue}
\usepackage{authblk}
\usepackage{url}
\usepackage{amsmath}
\usepackage{graphicx}
\usepackage{vmargin}
\setmarginsrb{2cm}{2cm}{2cm}{2cm}{0mm}{0mm}{0mm}{10mm}
\usepackage{enumitem}
\usepackage{booktabs}
\usepackage{marvosym}
\usepackage{tikz}
\usetikzlibrary{decorations.pathmorphing,patterns}
\usetikzlibrary{calc,patterns,decorations.pathmorphing,decorations.markings}

\title{\bfseries The effects of increasing speed on the structure\\ of axle box bogies of NT-11 trains}
\author[1]{\normalsize S. W. Rienstra}
\author[2]{N. Karjanto\thanks{Corresponding author \ \Letter: \protect\url{nkarjanto@gmail.com}.}}
\author[2]{S. D. Destiyani}
\author[2]{Andonowati}
\affil[1]{\small Department of Mathematics and Computer Science, Eindhoven University of Technology \authorcr Den Dolech 2, 5612 AZ, Eindhoven, The Netherlands}
\affil[2]{Department of Mathematics, Bandung Institute of Technology \authorcr Jalan Ganesha 10, Bandung 40132, Indonesia}
\date{}

\begin{document}
\maketitle

\begin{abstract}
The occurrence of cracks on the structure of axle box bogie of NT-11 economic trains operated by PT Kereta Api Indonesia is studied. These axle box bogies were designed for the maximum speed of 90~km/hour. The current market demand led to an increase in maximum speed to 110~km/hour. We investigate the effects of this speed increase in the occurrence of cracks of the axle box bogies. Based on the data obtained from PT Kereta Api Indonesia, a statistical approach is first used to see a relationship between the occurrence of cracks and some parameters, namely the distance to one-way destination, the speed, and the number of stops that the trains make to arrive at the destination. From this statistical study, stopping frequency is the most dominant factor that affects the occurrence of cracks and that the increase in speed leads to an increase in the crack's occurrence. A mathematical model is then designed to further describe the effects of this increased velocity on the structure of axle box bogies.  From this model, a resonance analysis is derived and numerical experiments are carried out. It is shown that the wagon resonance is highly dependent on the velocity. An example of typical data shows that this resonance occurs when the train reaches a velocity at $\approx 104$~km/hour. These resonance effects are believed to create a larger distance between parts of the axle box bogie that influence the occurrence of the cracks.
\end{abstract}

\section{Introduction}

The market increase in both passenger and cargo railways has led PT Kereta Api Indonesia as a company that manages the Indonesian railway business to consider an improvement in its capacity as well as its services. For certain corridors such as the dense railway traffics Jakarta--Bandung or Jakarta--Cirebon, to answer this market demand the increase in the number of trains should be balanced with the need for new tracks. On the other hand, there is a customer demand to reduce the time required for covering the travel distance by increasing the train's speed.

For the past few years, the passengers of the economy train NT-11 has increased quite drastically. The number of trains used, however, cannot balance the need for such carriers. This leads to an increase in the frequency of the train's operation as well as the train's speed.  This move is believed to cause some operational as well as structural problems. The occurrence of cracks on the axle box bogies leading to the broken structure of the NT-11 trains is an example of such problems.

In the beginning, the train axle box bogie of NT-11 was designed for the maximum speed of 90~km/hour. The current market demand for faster trains leads to an increase in speed up to 110~km/hour. Ever since this maximum operational speed was implemented, only in the past two years, there has been a failure of 60 NT-11 trains out of 200~trains. This article is an attempt to understand how the increase in speed can play a significant role in the occurrence of cracks of the axle box bogies of NT-11 trains.

One of the aims of such an understanding is to suggest to the PT Kereta Api Indonesia to reconsider reducing the need for high speeds that will serve to help in solving this crack problem. The high demand for such passenger carriers may be solved by redesigning the schedule of the train's operations~\cite{kai}, cf.~\cite{yos}.

This article is arranged in the following order. In the next section, we will present some figures describing the physical structures of the train's wagons together with their bogies. Some data collected from PT Kereta Api Indonesia on the occurrence of cracks and the related operational parameters are also described here. These parameters are the distance to the destination covered by the trains, the average operational speed, and the stopping frequency. Based on these collected data, a statistical approach is used in Section~\ref{statistic} to see the effects of these parameters on the occurrence of the cracks. It is shown here that the increase in the average operational speed influences the increase in the crack's percentage. In Section~\ref{model}, some simple mathematical models are designed and a resonance analysis is further studied. For some typical data, numerical experiments are carried out. It is shown that the operational speed highly influences the wagon resonance. This resonance is believed to cause the cracks leading to the broken structure and failure of the NT-11 trains. In the last section, we present some concluding remarks and suggestions.

\section{Some facts}			\label{fact}

\begin{figure}[htbp]
\begin{center}
\includegraphics[scale = 0.5]{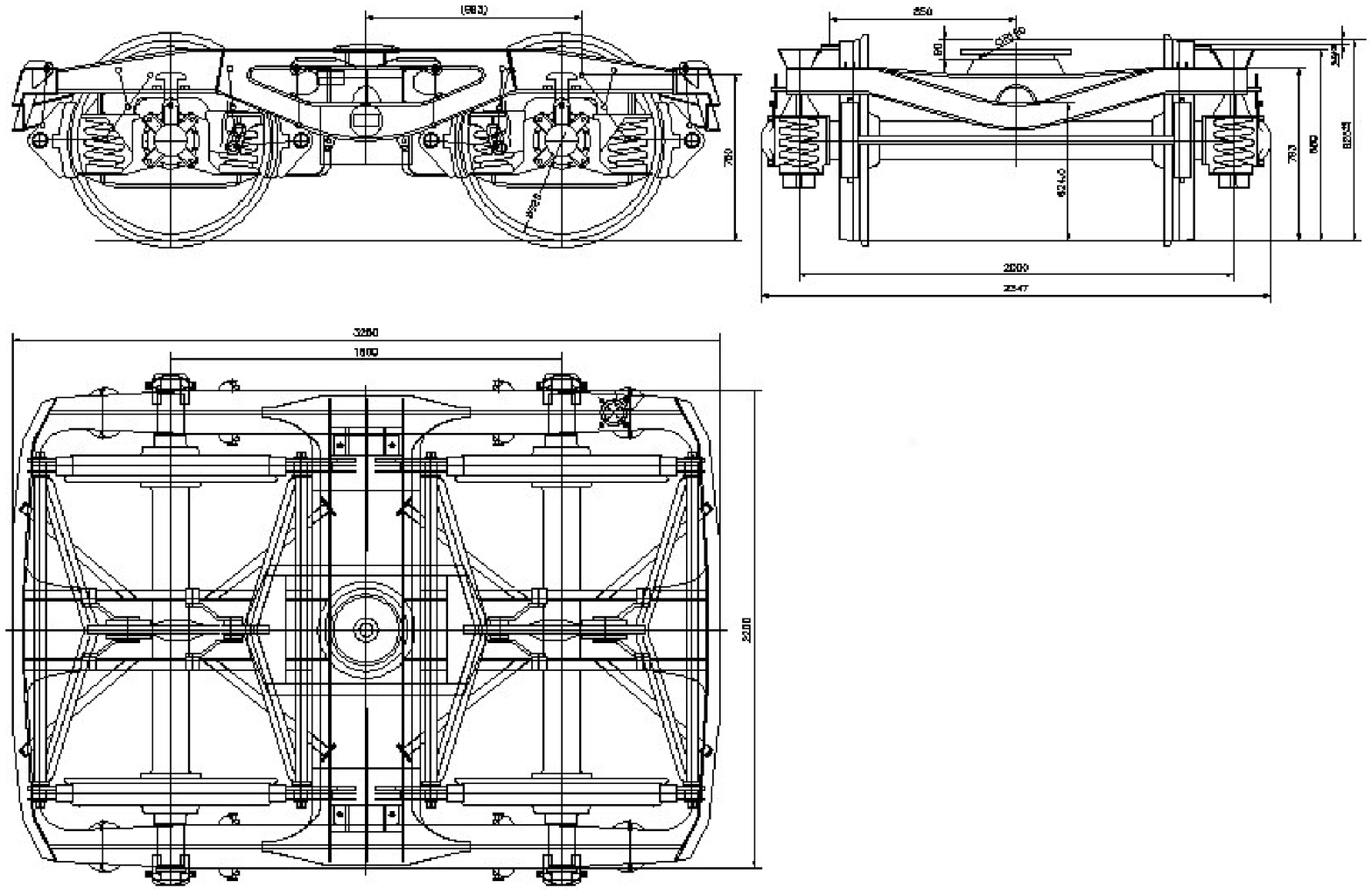}
\vspace*{-0.4cm}	
\end{center}
\caption{A structure of an axle box bogie consisting of a system of spring, a bent cylinder connecting the spring and a metal box covering the system. On the first row is the side view (left) and its details (right) of the axle bogie while the second row is the view of the bogie from above.} 				\label{structure}
\end{figure}

An axle box bogie is a part of a system in the railway's machinery that has an important role in the movement of traditional trains. Physically, this bogie consists of a supporting metal cylinder and a system of spring connecting the train's wagons and wheels. Every couple of spring is connected by a bogie frame where above this frame the train's wagon rests. The bogies of a train are located close to the center of the train's wheels and are covered by a box and connected by an axle (axle box bogie). This structure of the axle bogie box can be seen from Figure~\ref{structure}.

The economic trains operated by PT Kereta Api Indonesia have axle box bogie of type NT-11. This type was designed for the maximum speed of 90~km/hour. In current practice, the maximum speed has increased to 110~km/hour. Despite the maximum speed that will influence the acceleration as well as the deceleration, in general, it will change the force dynamics working on the trains.

In what follows, we present some data obtained for the economic trains operated in Java using the NT-11 axle box bogies. The data relate the percentage of the cracks in the axle box bogies, the one-way distance to the destination~(km), the number of stops and the average speed during which the trains are operated~(km/hour) in covering the distance. 
\begin{table}[h]
\begin{center}
\begin{tabular}{@{}ccccc@{}}
\toprule
$i^{\text{th}}$ observation & Percentage of cracks & Distance covered & Number of stops & Average speed \\ \hline
 1 & 0.140 & 54.81 &  9 &  95 \\
 2 & 0.350 & 79.69 & 16 &  90 \\
 3 & 0.125 & 71.14 &  3 & 110 \\
 4 & 0.125 & 47.16 &  5 & 100 \\
 5 & 0.125 & 43.86 &  4 & 100 \\
 6 & 0.130 & 51.22 &  6 &  95 \\
 7 & 0.125 & 63.14 &  8 &  90 \\
 8 & 0.125 & 46.16 &  6 & 100 \\
 9 & 0.200 & 69.31 & 10 &  95 \\
10 & 0.200 & 54.28 &  9 &  90 \\
11 & 0.302 & 61.29 &  9 &  95 \\
12 & 0.236 & 53.49 &  7 &  95 \\
13 & 0.260 & 49.86 &  5 & 100 \\
14 & 0.210 & 50.25 &  6 &  95 \\
15 & 0.165 & 49.06 &  5 & 110 \\
\bottomrule
\end{tabular}
\end{center}
\caption{An observation of the relationship between the percentage of the crack from the axle box bogies, the one-way distance to the destination (km), the number of stops and the average speed during which the trains are operated~(km/hour) to cover the distance.}   \label{tabcra}
\end{table}

\section{How increasing speed influences the percentage of cracks:\\ a statistical approach}     \label{statistic}

\subsection{A brief  description of linear regression model}       \label{linreg}

In this approach, we are going to use a linear regression model~\cite{reg}. The term linear here refers to the parameters in the model but not the free variables. For example, $Y = b_0 + b_1 x_1 + b_2 x_2^2 + b_3 \log x_3 + \epsilon$, for a dependent variable $Y$ and free variables $x_1$, $x_2$, and $x_3$, is a linear model in the parameters $b_0$, $b_1$, $b_2$, and $b_3$ (with an error $\epsilon$) although it is quadratic in $x_1$ and logarithmic in $x_3$. Thus the linear regression model comprises a large class of statistical models. In our case, the model will take the dependent variable percentage of cracks as a function of distance, stopping frequency and speed.

In principle, when the model is not linear in the free variables under consideration, we may define new free variables in such away that the new model is a linear function of the defined variables. In the above example, $Y = b_0 + b_1 x_1 + b_2 x_2^2 + b_3 \log x_3 + \epsilon$, we  may  take $z = x_1$, $z_2 = x_2^2$, and $z_3 = \log x_3$ to transform the model into into $Y = b_0 + b_1 z_1 + b_2 z_2 + b_3 z_3 + \epsilon$. So wee can consider a general linear regression model with $k$ free variables in the form $Y = \beta_0 + \beta_1 x_1 + \beta_2 x_2 + \dots + \beta_k x_k + \epsilon$.

Suppose that a sample consisting of $n$ observation describing a dependent variable $Y$ and $k$ free variables, $j = 1, 2, \dots, k$, are taken. Let $y_i$ and $x_{ij}$; $i = 1, 2, \dots, n$, $j = 1, 2, \dots, k$; denote the value $Y$ and $x_j$ ($j = 1, 2, \dots, k$) at the  observation. If $\epsilon_i$  is the error of this $i^{\text{th}}$ observation, then  substituting these data into the model we obtain a linear equation in the form:
\begin{equation}
\left[\begin{array}{c}
y_1 \\ y_2 \\ \vdots \\ y_n
\end{array}\right] = 
\left[\begin{array}{ccccc}
1 & x_{11} & x_{12} & \dots & x_{1k} \\
1 & x_{21} & x_{22} & \dots & x_{2k} \\
\vdots & \vdots & \vdots & \ddots & \vdots \\
1 & x_{n1} & x_{n2} & \dots & x_{nk}
\end{array}\right]
\left[\begin{array}{c}
\beta_0 \\ \beta_1 \\ \vdots \\ \beta_k
\end{array}\right] + 
\left[\begin{array}{c}
\epsilon_1 \\ \epsilon_2 \\ \vdots \\ \epsilon_n
\end{array}\right].   \label{linsysmat}
\end{equation}
The vector coefficient $\beta = \left[\beta_0 \quad \beta_1 \quad \dots \quad \beta_k \right]^T$ has to be obtained by minimizing the error using the least-squares method. If we write
\begin{equation*}
J = \sum_{i = 1}^{n} \epsilon_i^2 = \sum_{i = 1}^{n} \left(y_i - \beta_0 - \beta_1 x_{i1} - \beta_2 x_{i2} - \dots - \beta_k x_{ik} \right)^2
\end{equation*}
then we have to minimize $J$ with respect to the parameters $\beta_0$, $\beta_1$, $\beta_2$, \dots, $\beta_k$ which can be easily obtained by equating partial derivatives of $J$ with respect to $\beta_i$ to zero, that is
\begin{equation}
\begin{aligned}
\frac{\partial J}{\partial \beta_0} &= -2 \sum_{i = 1}^{n} \left(y_i - \beta_0 - \beta_1 x_{i1} - \beta_2 x_{i2} - \dots - \beta_k x_{ik} \right) = 0 \\
\frac{\partial J}{\partial \beta_1} &= -2 \sum_{i = 1}^{n} \left(y_i - \beta_0 - \beta_1 x_{i1} - \beta_2 x_{i2} - \dots - \beta_k x_{ik} \right) x_{i1} = 0 \\
\frac{\partial J}{\partial \beta_2} &= -2 \sum_{i = 1}^{n} \left(y_i - \beta_0 - \beta_1 x_{i1} - \beta_2 x_{i2} - \dots - \beta_k x_{ik} \right) x_{i2} = 0 \\
\vdots \quad & \qquad \qquad \qquad \qquad \qquad \qquad \qquad \vdots \\
\frac{\partial J}{\partial \beta_k} &= -2 \sum_{i = 1}^{n} \left(y_i - \beta_0 - \beta_1 x_{i1} - \beta_2 x_{i2} - \dots - \beta_k x_{ik} \right) x_{ik} = 0.
\end{aligned}  \label{parderJ}
\end{equation}
From~\eqref{parderJ}, it is not difficult to see that estimates for $\beta = \left[\,\beta_0 \quad \beta_1 \quad \dots \quad \beta_k \, \right]^T$ satisfy
\begin{equation*}
\left(\mathbf{X}^T \mathbf{X} \right) \beta = \mathbf{X}^T \mathbf{Y}
\end{equation*}
with
\begin{equation*}
\mathbf{X} = \left[\begin{array}{ccccc}
1 & x_{11} & x_{12} & \dots & x_{1k} \\
1 & x_{21} & x_{22} & \dots & x_{2k} \\
\vdots & \vdots & \vdots & \ddots & \vdots \\
1 & x_{n1} & x_{n2} & \dots & x_{nk}
\end{array}\right] \qquad \text{and} \qquad
\mathbf{Y} = \left[\, y_1 \quad y_2 \quad \dots \quad y_n \, \right]^T.
\end{equation*}
If $\mathbf{X}^T \mathbf{X}$ is not singular, so that the inverse exists, then $\beta$ can be calculated as
\begin{equation}
\beta = \left(\mathbf{X}^T \mathbf{X}\right)^{-1} \mathbf{X}^T \mathbf{Y}.  \label{invbet}
\end{equation}

\subsection{Application of linear regression to the crack problem}   \label{apcra}

In this application, the dependent variable is the percentage of crack axle bogie box for NT-11 type. Let \texttt{MinDist}, \texttt{MinStop}, and \texttt{MinAverSpeed} be the values (obtained from data) of the minimum distance, the minimum number of stops, and the minimum average speed,  respectively.  We would like to see the effect of the increasing distance, the number of stops, and the speed on the percentage of cracks. Here as  the free variables, we take the differences between the distance and \texttt{MinDist}, stopping frequency and \texttt{MinStop}, the average speed, and \texttt{MinAverSpeed}.

Let us consider a model $Y = \beta_0 + \beta_1 x_1 + \beta_2 x_2 + \beta_3 x_3 + \epsilon$. Here $Y$ represents the percentage of cracks, while $X_j$, $j = 1, 2, 3$ are the differences between the distance and \texttt{MinDist}, stopping frequency and \texttt{MinStop}, the average speed and \texttt{MinAverSpeed}, respectively. The sample in Table~\ref{tabcra} consisting of 15 observations describing the dependent variable $Y$ and three free variables $x_j$, $j = 1, 2, 3$. Let $y_i$ and $x_{i}$; $i = 1, 2, \dots, 15$, $j = 1, 2,3$; denote the value $Y$ and $x_j$ ($j = 1, 2, 3$) at the $i^{\text{th}}$ observation. If $\epsilon_i$ is the error of this $i^{\text{th}}$ observation, then matrix $\mathbf{X}$ and the vector $\mathbf{Y}$ in the linear system of equations~\eqref{linsysmat} are given by
\begin{equation*}
\mathbf{X} = \left[\begin{array}{cccc}
1 & x_{1,1} & x_{1,2} & x_{1,3} \\
1 & x_{2,1} & x_{2,2} & x_{2,3} \\
\vdots & \vdots & \vdots & \vdots \\
1 & x_{15,1} & x_{15,2} & x_{15,3}
\end{array}\right] \qquad \text{and} \qquad
\mathbf{Y} = \left[\, y_1 \quad y_2 \quad \dots \quad y_{15} \, \right]^T.
\end{equation*}
where the second, third and fourth columns of $\mathbf{X}$ are obtained from the columns of ``Distance covered'', the ``Number of stops'', and the ``Average speed'' in Table~\ref{tabcra}, after subtraction by the minimum distance, the minimum number of stops, and minimum average speed, respectively.  The vector $\mathbf{Y}$ is obtained from the column of the ``Percentage of cracks'' in Table~\ref{tabcra}. The straightforward application of~\eqref{invbet} gives
\begin{equation*}
\beta 
= \left[\, \beta_0 \quad \beta_1 \quad \beta_2 \quad \beta_3 \, \right]^T 
= \left(\mathbf{X}^T \mathbf{X}\right)^{-1} \mathbf{X}^T \mathbf{Y} 
= \left(0.097, \; -0.049, \; 0.018, \; 0.0021 \right).
\end{equation*}
Thus,
\begin{equation}
Y = 0.097 - 0.049 X_1 - 0.018 X_2 + 0.0021 X_3 + \text{error}  \label{regval}
\end{equation}
where the relative error of this computation is 
\begin{equation*}
\text{Relative Error} = \frac{1}{\sum_{i = 1}^{n} \left(y_i\right)^2} \sum_{i = 1}^{n} \left(y_i - \beta_0 - \beta_1 x_{i1} - \beta_2 x_{i2} - \dots - \beta_k x_{ik} \right)^2 \approx 6.6\%.
\end{equation*}
From~\eqref{regval}, it shows here that the one-way distance covered produces a negative but rather insignificant effect on the percentage of cracks. The longer the distance covered, the less likely the cracks occurred on the axle box bogie. The number of stops during the operation of the trains is a positive and gives the most significant effect on the occurrence of the cracks. If stopping frequency is defined as the number of stops divided by the distance covered, it can be concluded that the percentage of cracks increases as the stopping frequency increases. Finally, from this observation, it shows that the increase in train speed increases the percentage of cracks.

\section{Some mathematical models}  \label{model}

\subsection{Train motion and its scaling rules}  \label{motion}

The profile of a railway depends on the geographical condition of its location. Between Bandung and Jakarta for example, the Bandung--Padalarang path is relatively flat while between Padalarang--Cikampek where it is a mountainous area, the railway profile is affected by this condition. The motion of a train follows from the railway where it operates. Since the position of the rails between Bandung and Jakarta varies slightly along the track, this train motion is not just one of horizontal motion, rather there is some vertical motion according to the profile of the rails.
\begin{figure}[htbp]
\begin{center}
\vspace*{-1cm}
\begin{tikzpicture}
\draw[thick] (0,0) .. controls (5,3) and (6,-3) .. (14,0.5);
\draw[blue,<->,dashed] (4,0.8)--(14,0.8);
\draw[blue] (9,1) node[above]{$L$};
\draw[red,<->,dashed] (9.25,0.75)--(9.25,-0.6);
\draw[red](9.25,0) node[left]{$A$};
\draw[dashed,->] (3,-0.2) .. controls (3.5,0.5) and (4,0) .. (4.8,0.4);
\draw (3,0) node[below]{${\displaystyle y = A h\left(\frac{x}{L}\right)}$};
\end{tikzpicture}
\vspace*{-2cm}
\end{center}
\caption{The profile of a segment of a railway.}	\label{railway}
\end{figure}
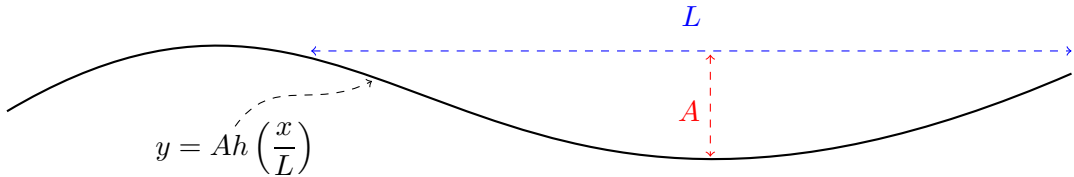

We consider a segment of a railway having a profile with some characteristic length~$L$ and amplitude~$A$ representing the vertical variation of the segment. (See Figure~\ref{railway}.) Let us consider a point on the train that is moving along this segment. If $u$ is the (fixed) horizontal speed of the train, then the position vector of this point can be written in the form
\begin{equation}
\mathbf{r}(t) 
= x(t) \, \mathbf{i} + y(t) \, \mathbf{j} 
= u t \, \mathbf{i} + A h\left(\frac{x}{L}\right) \, \mathbf{j} 
= u t \, \mathbf{i} + A h\left(\frac{u t}{L}\right) \, \mathbf{j}.  \label{posivec}
\end{equation}
Here, $h$ is a function describing the profile of the segment of the rails and $t$ is time.

From~\eqref{posivec}, the velocity and acceleration of the point while passing the segment are 
\begin{subequations}
\renewcommand{\theequation}{\theparentequation.\arabic{equation}}
\begin{align}
\mathbf{v}(t) & = \dot{\mathbf{r}}(t) = \dot{x}(t) \, \mathbf{i} + \dot{y}(t) \, \mathbf{j} 
				= u \, \mathbf{i} + A \frac{u}{L} h'\left(\frac{u t}{H}\right) \, \mathbf{j}\\
\mathbf{a}(t) & = \dot{\mathbf{v}}(t) = \ddot{\mathbf{r}}(t) = \ddot{x}(t) \, \mathbf{i} + \ddot{y}(t) \, \mathbf{j}
                = A \frac{u^2}{L^2} h''\left(\frac{u t}{H}\right) \, \mathbf{j}.  \label{accele}
\end{align}
\end{subequations}
The final expression of~\eqref{accele} shows that the vertical acceleration the point under consideration scales on the square of the horizontal speed of the train, i.e. $\ddot{y} = {\mathcal{O}}\left(A \frac{u^2}{L^2}\right)$.  Although the horizontal acceleration is zero, this vertical component influences the forces working on the bogie system as it will be shown in the following paragraph.
\begin{figure}[htbp]
\begin{center}
\begin{tikzpicture}
\draw[thick,rounded corners,red,fill = red!20] (0, 0) rectangle (0.3, 2) {};
\draw[very thick] (0.3, 1) -- (5.7, 1) {};
\draw[thick,rounded corners,red,fill = red!20] (5.7, 0) rectangle (6, 2) {};
\draw[<->, dashed] (0.4,0.1) -- (5.6,0.1);
\draw (3,0) node[below]{$B$};
\draw[->, dashed] (0.15,-1) -- (0.15,-0.1); 
\draw (0.2,-1) node[below]{$y_1$};
\draw[->, dashed] (5.85,-1) -- (5.85,-0.1); 
\draw (5.9,-1) node[below]{$y_2$};
\draw[very thick,blue, fill = blue!30] (3,3) circle [radius = 0.2];
\draw[<->,dashed] (3,1.1) -- (3,2.75);
\draw (3,2) node[right]{$H$};
\draw[<->, dashed] (2.1,3.25) -- (3.9,3.25);
\draw (3,3.25) node[above]{$w$};
\draw[<->,thick, dashed] (6.5,0.5) to [bend right = 40] (6.5,1.5);
\draw (6.7,1) node[right]{$\alpha$};
\end{tikzpicture}
\vspace*{-0.5cm}
\end{center}	
\caption{A simplification of  the structure of Figure~\ref{structure}. The train's wheels are separated by a distance of~$B$ units, the wagon is located~$H$ units above the wheels and gives pressure to the system of axle bogie located in between the wheels. The blue circle is a point of reference for the position of the wagon. Here, $\alpha$~is the angle between the wheels' position while the train is on the segment under consideration and the horizontal position. Additionally, $y_1$~and $y_2$~are the vertical ordinate of the wheels' positions.}	 \label{simplif}
\end{figure}
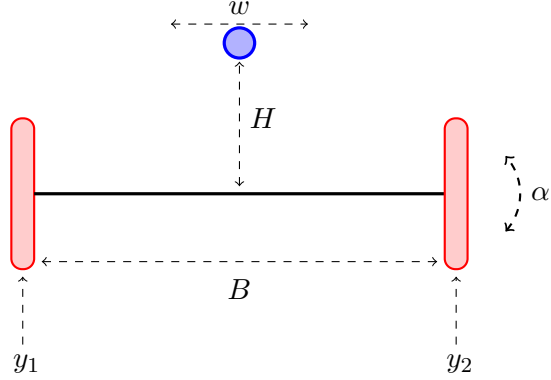

We consider a simplified structure of a system consisting of two train's wheels, a bogie, and a wagon (See Figures~\ref{structure} and~\ref{simplif}). Suppose that the wheels are separated by a distance~$B$ and the wagon is located~$H$ units above the wheels. A system of an axle bogie is located in between the wheels. On a flat rail, this system is pictured in Figure~\ref{simplif}. When the train moves to the segment as shown in Figure~\ref{railway}, the wheels' position makes an angle, say~$\alpha$, from its horizontal position. Suppose that $w$~is the position vector of a point in the wagon with respect to the point of reference as shown in Figure~\ref{simplif}. Then as the train enters the rail segment, it results in the change of the wheels' positions. This change further affects the change in $w$ as follows:
\begin{equation*}
\tan \alpha = \frac{y_2 - y_1}{B} = \frac{A}{B} \tilde{h} \left(\frac{ut}{\tilde{L}}\right)
\end{equation*}
for some function $\tilde{h}$ related to the function $h$ as described in~\eqref{posivec}. The change in $w$ can be written as
\begin{equation*}
w = H \tan \alpha = H \frac{A}{B} \tilde{h} \left(\frac{ut}{\tilde{L}}\right)
\end{equation*}
resulting to
\begin{equation}
\ddot{w} = {\mathcal{O}}\left(H \frac{A}{B} \frac{u^2}{\tilde{L}^2} \right).
\end{equation}

\subsection{Resonance analysis on a model of a spring-mass system}  \label{resonanceanalysis}

In this subsection, we will investigate the effects of the increase in the train's speed on the cracks of the axle box bogie as a resonant phenomenon. For this reason, we model the problem as a mass-spring system consisting of the wagon and the bogie~\cite{kreyszig}. We assume that the mass of the wheel is relatively negligible. Because only the wagon and the bogie are included, we call this model a simple two-mass model.
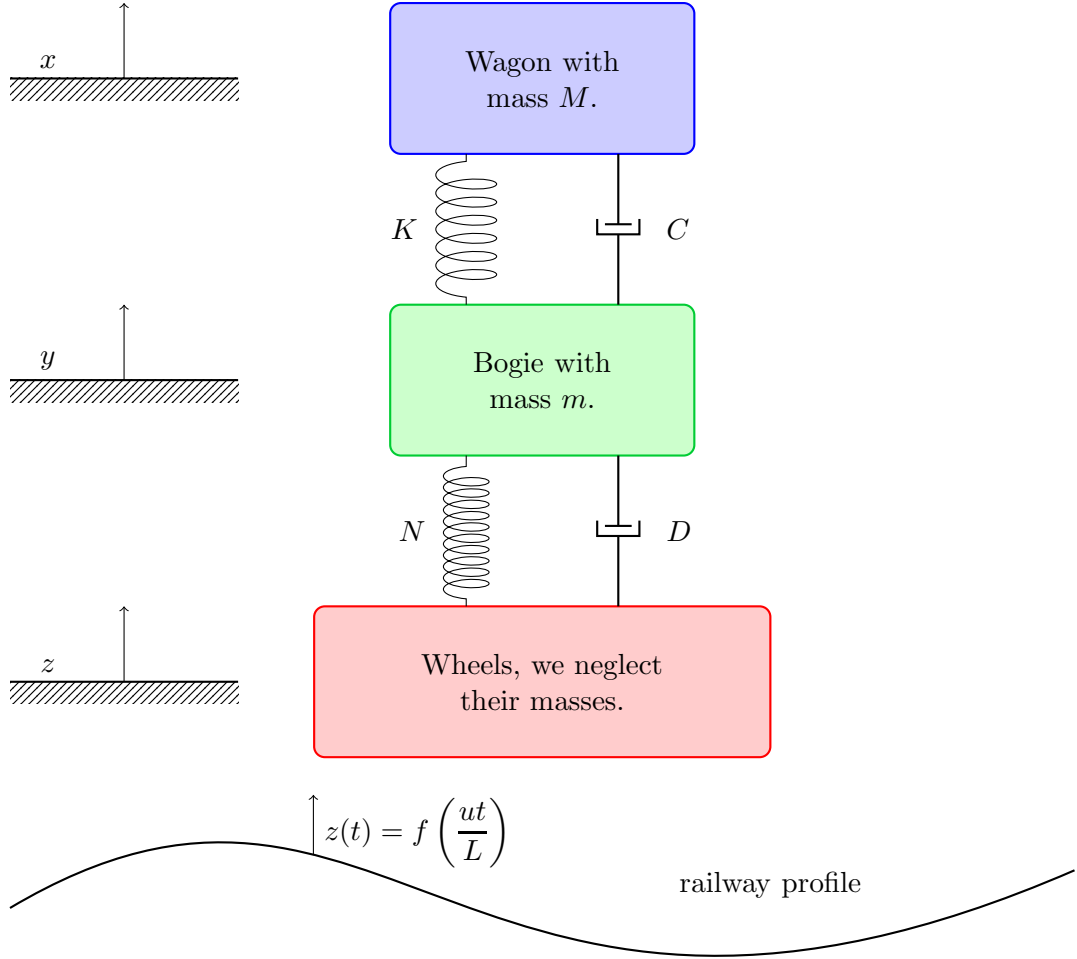
\begin{figure}[htbp]
\begin{center}
\begin{tikzpicture}
\draw[thick] (0,0) .. controls (5,3) and (6,-3) .. (14,0.5);
\draw (10,0) node[above]{railway profile};
\draw[->] (3.99,0.7) -- (3.99,1.5);
\draw (4,1) node[right]{${\displaystyle z(t) = f\left(\frac{ut}{L}\right)}$};
\draw[thick,rounded corners,red,fill = red!20,text = black, text width = 3.5cm, align = center] (4,2) rectangle node{Wheels, we neglect their masses.} (10,4); 
\draw[thick,rounded corners,green!80!blue, fill = green!20, text = black, text width = 2cm, align = center] (5,6) rectangle node{Bogie with mass $m$.} (9,8);
\draw[thick,rounded corners,blue, fill = blue!20, text = black, text width = 2.1cm, align = center] (5,10) rectangle node{Wagon with mass $M$.} (9,12);
\draw[decoration={aspect=0.3, segment length=1.5mm, amplitude=3mm, coil, pre length = 1mm, post length = 1mm},decorate] (6,4) -- (6,6);
\draw (5.6,5) node[left]{$N$};
\draw[decoration={aspect=0.3, segment length=2.4mm, amplitude=4mm, coil, pre length = 1mm, post length = 1mm},decorate] (6,8) -- (6,10);
\draw (5.5,9) node[left]{$K$};
\tikzstyle{damper}=[thick,decoration={markings,  
	mark connection node=dmp,
	mark=at position 0.5 with 
	{
		\node (dmp) [thick,inner sep=0pt,transform shape,rotate=-90,minimum width=15pt,minimum height=3pt,draw=none] {};
		\draw [thick] ($(dmp.north east)+(2pt,0)$) -- (dmp.south east) -- (dmp.south west) -- ($(dmp.north west)+(2pt,0)$);
		\draw [thick] ($(dmp.north)+(0,-5pt)$) -- ($(dmp.north)+(0,5pt)$);
	}
}, decorate]
\draw[damper] (8,4) -- (8,6);
\draw (8.5,5) node[right]{$D$};
\draw[damper] (8,8) -- (8,10);
\draw (8.5,9) node[right]{$C$};
\tikzstyle{ground}=[fill,pattern=north east lines,draw=none,minimum width=3cm,minimum height=0.3cm,xshift=0.5cm]
\node [ground,anchor=south] at (1,2.7){};
\draw[thick] (0,3) -- (3,3);
\draw[->] (1.5,3) -- (1.5,4);
\draw (0.5,3) node[above]{$z$};
\node [ground,anchor=south] at (1,6.7){};
\draw[thick] (0,7) -- (3,7);
\draw[->] (1.5,7) -- (1.5,8);
\draw (0.5,7) node[above]{$y$};
\node [ground,anchor=south] at (1,10.7){};
\draw[thick] (0,11) -- (3,11);
\draw[->] (1.5,11) -- (1.5,12);
\draw (0.5,11) node[above]{$x$};
\end{tikzpicture}
\vspace*{-2.5cm}
\end{center}
\caption{A spring-mass system consisting of the train's wagon and bogie. The mass of the wheels is neglected. Several data for the units in this figure are given $M = 11.2$~tons, $m = 1.01$~tons, $K = 963$~kN/m, $N = 1540$~kN/m, $C = 50$~kNs/m, $D = 0$.}
\end{figure}

The mathematical model for the two-mass system is given by
\begin{equation}
\left\{
\begin{aligned}
M \ddot{x} + K(x - y) + C(\dot{x} - \dot{y}) &= 0 \\
m \ddot{y} + N(y - z) + D(\dot{y} - \dot{z}) - K(x - y) - C(\dot{x} - \dot{y}) &=0.
\end{aligned}
\right. \label{twomass}
\end{equation}
As the model is linear, it is no restriction to consider the geometry perturbation function~$f$ in its Fourier components. Let us consider the dominating wavelength, which will for any reasonably produced railway track, typically be of the order of the length~$L$ of a section. This is equivalent to saying that (for this Fourier component)  the railway profile is expressed in term of  the sinusoidal function of period one. Writing $f(x) = A \cos(2\pi x) = \text{Re} \left\{A \cos(2\pi x) + i \sin (2\pi x) \right\} = \text{Re}\left\{A e^{2 i \pi x} \right\}$, then ${\displaystyle f\left(\frac{ut}{L}\right) = \text{Re} \left\{A e^{i \omega t} \right\}}$, where ${\displaystyle \omega = \frac{2\pi}{L} u}$.

As the problem is linear, causality arguments tell us that the response (the  solution of the system of  differential equations~\eqref{twomass}) has the same frequency~$\omega$, i.e., of the form
\begin{equation*}
\left\{
\begin{aligned}
x &= x_0 e^{i \omega t} \\
y &= y_0 e^{i \omega t} \\
z &= A e^{i \omega t}.
\end{aligned}
\right.
\end{equation*}
Substituting these expressions to the differential equations system~\eqref{twomass}, we obtain
\begin{equation*}
\left\{
\begin{aligned}
-\omega^2 M x_0 + (K + i \omega C) (x_0 - y_0) &= 0 \\
-\omega^2 m y_0 + (N + i \omega D) (y_0 - A) - (K + i \omega C) (x_0 - y_0) &= 0.
\end{aligned}
\right.
\end{equation*}
Rewriting these equations in the matrix form, we have
\begin{equation}
\left[
\begin{array}{cc}
-\omega^2 M + (K + i \omega C) & -(K + i \omega C)  \\
             -(K + i \omega C) & -\omega^2 m + (N + i \omega D) + (K + i \omega C)  
\end{array}
\right]
\left[
\begin{array}{c}
x_0 \\ y_0
\end{array}
\right] = 
\left[
\begin{array}{c}
0 \\ (N + i\omega D) A
\end{array}
\right]. 			\label{matrixsystem}
\end{equation}
Let us use the symbols $\tilde{K} = K + i \omega C$ and $\tilde{N} = N + i \omega D$ in the matrix expression above. Using these symbols, the matrix equation~\eqref{matrixsystem} simplifies to 
\begin{equation*}
\underbrace{\left[
\begin{array}{cc}
\omega^2 M - \tilde{K} & \tilde{K}  \\
\tilde{K} & \omega^2 m - \tilde{N} - \tilde{K}  
\end{array}
\right]}_S
\left[
\begin{array}{c}
x_0 \\ y_0
\end{array}
\right] = 
\left[
\begin{array}{c}
0 \\ -\tilde{N} A
\end{array}
\right].
\end{equation*}
Consequently, the solution of the matrix equation is given by
\begin{equation}
\left[
\begin{array}{c}
x_0 \\ y_0
\end{array}
\right] = 
\frac{1}{\text{det}(S)}
\left[
\begin{array}{cc}
\omega^2 m - \tilde{N} - \tilde{K} & -\tilde{K}  \\
-\tilde{K} & \omega^2 M - \tilde{K}
\end{array}
\right]
\left[
\begin{array}{c}
0 \\ -\tilde{N} A
\end{array}
\right]				\label{solution}
\end{equation}
where det$(S) = (\omega^2 M - \tilde{K})(\omega^2 m - \tilde{N} - \tilde{K}) - \tilde{K}^2$.

If det$(S)$ in~\eqref{solution} tends to zero, then the displacements $x_0$ and $y_0$ tend to infinity. Physically, the wagon and the bogie shake extremely causing the springs connecting them to break. This condition is called resonance. To simplify our computation, we observe some simple cases, such as for $C = D = 0$ and $m/M = \epsilon \ll 1$. The resonance will happen if det$(S) = 0$, i.e.,
\begin{equation}
\text{det}(S) = (\omega^2 M - K)(\omega^2 m - N - K) - K^2 = 0.  \label{resonance}
\end{equation}
Introducing new notations
\begin{equation*}
\Omega^2 = \omega^2 \frac{M}{K} \qquad \qquad \text{and} \qquad \qquad
\Lambda = \frac{N + K}{K}
\end{equation*}
equation~\eqref{resonance} above becomes
\begin{equation*}
K^2 \left(\omega^2 \frac{M}{K} - 1 \right) \left(\frac{m}{M} \omega^2 \frac{M}{K} - \frac{N + K}{K} \right) - K^2 = 0
\end{equation*}
resulting in
\begin{equation}
K^2 \left[ \left(\Omega^2 - 1\right) \left(\epsilon \Omega^2 - \Lambda \right) - 1 \right] 
\quad \qquad \Longrightarrow \qquad \quad
\left(\Omega^2 - 1\right) \left(\epsilon \Omega^2 - \Lambda \right) - 1 = 0.  			\label{key}
\end{equation}
\begin{itemize}[leftmargin=1.2em]
\item Case I: $\Omega^2 = {\mathcal{O}}(1)$
\begin{equation*}
\Omega^2 = 1 - \frac{1}{\Lambda} \quad \Longrightarrow \quad
\omega^2 \frac{M}{K} = 1 - \frac{K}{N + K} \quad \Longrightarrow \quad
\omega^2 = \frac{1}{M} \frac{NK}{N + K} \quad \Longrightarrow \quad
f_L = \frac{1}{2\pi} \sqrt{\frac{1}{M} \frac{NK}{N + K}}.
\end{equation*}
Substituting $M = 11.2$~tons; $N = 1540$~kN/m; and $K = 963$~kN/m into (12), the resonance frequency is $1.158$~Hz. If one rail section is $L = 25$~m, then the speed causing resonance is $u = f L = 1.158 \times 25 = 28.95$~m/sec $= 104.22$~km/h. This value will be passed when the train changes its velocity from 90~km/h to 110~km/h, so the tendency is that resonance is very well possible. This low-frequent resonance is a perturbation of the resonance of the heaviest element in the system, i.e., the wagon.

\item Case II: $\Omega^2 = {\mathcal{O}}\left(1/\epsilon\right)$
\begin{equation*}
\Omega^2 = \frac{\Lambda}{\epsilon} \quad \Longrightarrow \quad
\omega^2 \frac{M}{K} = \frac{N + K}{K} \, \frac{M}{m} \quad \Longrightarrow \quad
\omega^2 = \frac{N + K}{m} \quad \Longrightarrow \quad
f_h = \frac{1}{2\pi} \sqrt{\frac{N + K}{m}}.
\end{equation*}
Substituting $N = 1540$~kN/m; $K = 963$~kN/m; and $m = 1.01$~tons, so the resonance frequency is $7.923$~Hz. If the length of one rail section is $L = 25$~m, then the velocity causes resonance is $u = f L = 7.923 \times 25 = 198.075$~m/sec $= 713.07$~km/h. This is related to the resonance of the lightest element, i.e., the bogie. As the corresponding velocity is extremely high, this resonance is of no concern here.
\end{itemize}

\section{Concluding remarks} 			\label{conclusion}

We have investigated the crack problem in the axle box of NT-11 bogies statistically by employing the linear regression model. From the measured observation, the percentage of crack occurrence is significantly influenced by the stopping frequency during the train operation as well as the increase in the train's speed. However, the covered distance gives a negative and insignificant effect on the percentage of cracks.

An elementary mass-spring-damper system is analyzed of two masses excited by the forcing due to a typical motion along the rails. The two masses differ such that a low and a high resonance frequency are found. The lower one corresponds to a typical train velocity of about 100~km/h, which is remarkably well in agreement with the practical observations. This supports our conclusion that the rapidly increasing wear of material and the occurrence of cracks, when the train operating speed is increased above design conditions, are mainly due to the resonance of the wagon-bogie system, excited by the non-flat rail bed.

\section*{\large Acknowledgment}

{\small The authors gratefully acknowledge Akhmad Zikri (PT Kereta Api Indonesia, the Indonesian Railway Company), Edy Soewono, Hadi Susanto, Anna Tania Siska, Andres Suryadi, Sondang Marsinta Uli Panggabean, Lylye Sulaeman Yusup (Bandung Institute of Technology), Barbera Wilhelmina van de Fliert (University of Twente), Neville Fowkes (The University of Western Australia), and Wikaria Gazali (Bina Nusantara University) for the advice, suggestion and fruitful discussion not only during the Industrial Mathematics Week (held at Bandung Institute of Technology on 10--14 July 2000) but also afterward. \par}

\end{document}